\documentclass[number,citesort,dvips]{arxbj}
\usepackage{cursive}  

\aid{0}
\volume{16}
\issue{4}
\pubyear{2010}
\firstpage{1343}
\lastpage{1368}
\doi{10.3150/09-BEJ247}

\makeatletter
\def\eqref#1{(\ref{#1})}
\newtheorem{teorem}{Theorem}[section]
\newtheorem{lem}[teorem]{Lemma}
\newtheorem{prop}[teorem]{Proposition}
\newremark{rem}{Remark}[section]
\makeatother

\begin{document}
\begin{frontmatter}

\title{Criteria for hitting probabilities
with applications to systems of stochastic wave~equations}
\runtitle{Criteria for hitting probabilities}

\begin{aug}
\author[1]{\fnms{Robert C.} \snm{Dalang}\thanksref{1}\ead[label=e1]{robert.dalang@epfl.ch}\corref{}} \and
\author[2]{\fnms{Marta} \snm{Sanz-Sol\'e}\thanksref{2}\ead[label=e2]{marta.sanz@ub.edu}}
\runauthor{R.C. Dalang and M. Sanz-Sol\'e}
\address[1]{Institut de Math\'ematiques, Ecole Polytechnique F\'ed\'erale,
Station 8, 1015 Lausanne, Switzerland. \mbox{\printead{e1}}}
\address[2]{Facultat de Matem\`atiques, Universitat de Barcelona,
Gran Via 585, 08007 Barcelona, Spain.\\ \mbox{\printead{e2}}}
\pdfauthor{Robert C. Dalang, Marta Sanz-Sole}
\end{aug}

\received{\smonth{2} \syear{2009}}
\revised{\smonth{12} \syear{2009}}

%
\begin{abstract}
We develop several results on hitting probabilities of random fields
which highlight the role of the dimension of the parameter space. This
yields upper
and lower bounds in terms of Hausdorff measure and Bessel--Riesz
capacity, respectively.
We apply these results to a system of stochastic wave equations in
spatial dimension $k\ge1$
driven by a $d$-dimensional spatially homogeneous additive Gaussian
noise that is white in
time and colored in space.
\end{abstract}

%

\end{frontmatter}
%

%
\section{Introduction}
\label{s1}

There have recently been several papers on hitting probabilities for
systems of stochastic partial differential equations (SPDEs). The first
seems to be \cite{mt}, which mainly studied polarity of points for the
Gaussian random field which is the solution of a system of linear heat
equations in spatial dimension one, driven by space--time white noise.
Next, the paper \cite{dn04} studied hitting probabilities for a
nonlinear system of (reduced) wave equations in spatial dimension one
and established upper and lower bounds on hitting probabilities in
terms of Bessel--Riesz capacity.

The paper \cite{dkn07} considered a system of nonlinear heat equations
in spatial dimension $k=1$ with additive space--time white noise and
established lower and upper bounds on the probability that the solution
$(u(t,x), (t,x) \in\mathbb{R}_+ \times[0,1])$ hits a set $A \subset
\mathbb{R}
^d$ in terms of capacity and Hausdorff measure, respectively. In \cite
{dkn09}, these results were extended to systems of the same heat
equations, but with multiplicative noise. The paper \cite{dkn10}
extends these results to systems of nonlinear heat equations in spatial
dimensions $k\geq1$, driven by spatially homogeneous noise that is
white in time. Some other results on hitting probabilities for
parabolic SPDEs with reflection are contained in the papers \cite
{z02,z03,dmz06}.

The objective of this paper is to begin a similar program for systems
of stochastic wave equations, starting with the analog of \cite
{dkn07}. We note that properties of solutions of stochastic wave
equations in spatial dimensions $k >1$ are often much more difficult to
obtain than their analogs for heat equations, due to the greater
irregularity of the fundamental solution of the wave equation. One
example of this is the study in \cite{dss} of H\"older continuity of
sample paths for the $3$-dimensional wave equation.

In \cite{dkn07}, various conditions on the density of the random
vector $(u(t,x),  u(s,y))$ were identified that imply upper and lower
bounds on hitting probabilities. The conditions were expressed using a
``parabolic metric'' and were designed to be applied to the stochastic
heat equation driven by space--time white noise. They were applied there
first to study the linear stochastic heat equation; the nonlinear
stochastic heat equation with additive noise was then handled by
appealing to Girsanov's theorem. Because of the absence of a suitable
Girsanov's theorem for heat or wave equations in spatial dimensions
$k>1$ (a problem also noted in \cite{dkn10}), we will first develop
some general results that will also be useful for nonlinear equations.
In contrast with \cite{dkn07}, these results are designed to be used
for stochastic wave equations. We will apply them to linear wave
equations in spatial dimension $k\geq1$, driven by spatially
homogeneous noise that is white in time. In a forthcoming work, we
intend to use these general results to study the nonlinear stochastic
wave equation with additive and/or multiplicative noise.

More precisely, we consider here the $d$-dimensional stochastic process
$U=\{(u_i(t,x), i=1,\ldots,d), (t,x)\in[0,T]\times\mathbb{R}^k\}$ which
solves the system of SPDEs
\begin{equation}
\label{1.1}
\frac{\partial^2 u_i}{\partial t^2 }(t,x) - \Delta u_i(t,x)= \sum
_{j=1}^d \sigma_{i,j} \dot F^j(t,x)
\end{equation}
for $(t,x)\in  \,]0,T]\times\mathbb{R}^k$, with initial conditions
\begin{equation}
\label{1.2}
u_i(0,x)=\frac{\partial u_i}{\partial t}(0,x)=0.
\end{equation}
Here, $\Delta$ denotes the Laplacian on $\mathbb{R}^k$ and $\sigma
=(\sigma
_{i,j})$ is a deterministic, invertible, $d\times d$ matrix. The noise
process $\dot F:=(\dot F^1,\ldots,\dot F^d)$ is a centered
(generalized) Gaussian process whose covariance is informally given by
an expression such as
\begin{equation}\label{e1}
E(\dot F^i(t,x)\dot F^j(s,y)) = \delta_{i,j}  \delta(t-s)  \Vert
x-y\Vert^{-\beta},
\end{equation}
where $\delta_{i,j}$ denotes the Kronecker symbol, $\delta(\cdot)$
is the Dirac delta function at zero and $\beta>0$.
More precisely, let $C^\infty_0(\mathbb{R}^{k+1})$ denote the space of
infinitely differentiable functions with compact support and consider a
family of centered Gaussian random vectors $F=(F(\varphi) =
(F^1(\varphi),\dots,F^d(\varphi)),  \varphi\in C^\infty_0(\mathbb{R}
^{k+1}))$, with covariance function
\begin{equation}\label{e2}
E(F(\varphi)F(\psi))= \int_{\mathbb{R}_+} \mathrm{d}r \int_{\mathbb{R}^k}
\Gamma(\mathrm{d}x)
\bigl(\varphi(t,\cdot) \ast\tilde\psi(t,\cdot)\bigr)(x),
\end{equation}
where $\tilde\psi(t,x) := \psi(t,-x)$ and $\Gamma$ is a
non-negative and non-negative definite tempered measure on $\mathbb
{R}^k$. We
note that \eqref{e2} reduces to \eqref{e1} if $\Gamma(\mathrm{d}x) = \Vert
x\Vert^{-\beta} \,\mathrm{d}x$. By the Bochner--Schwartz theorem (see \cite
{s}), there exists a non-negative tempered measure $\mu$ on $\mathbb{R}^k$
(termed the \textit{spectral measure} of $F$) such that $\Gamma= \mathcal
{F}\mu
$, where $\mathcal{F}$ denotes the Fourier transform. Elementary
properties of
the Fourier transform show that equation \eqref{e2} can be written
\begin{equation}
\label{e20}
E(F(\varphi)F(\psi))= \int_{\mathbb{R}_+} \mathrm{d}r \int_{\mathbb{R}^k}
\mu(\mathrm{d}\xi)
\mathcal{F}\varphi(t,\cdot)(\xi)  \overline{\mathcal{F}\psi
(t,\cdot)(\xi)}.
\end{equation}
Let $G(t,x)$ be the fundamental solution of the wave equation.
Generically, the solution $u$ of \eqref{1.1} is given by
\begin{equation}
\label{1.3}
u_i(t,x)=\int_0^t\int_{\mathbb{R}^k}G(t-r,x-y) \sum_{j=1}^d \sigma_{i,j}
M^j(\mathrm{d}r,\mathrm{d}y) ,
\end{equation}
where $M=(M^1,\ldots,M^d)$ is the martingale measure derived from
$\dot F$ (see \cite{df} for details). However, it is well known that
$G$ is a function in dimensions $k\in\{1,2\}$ only, so the stochastic
integral in \eqref{1.3} should be interpreted in the sense of \cite
{d99}. We note that, according to (\ref{e2}) and (\ref{e20}),
\[
E((u_i(t,x))^2)= \Biggl(\sum_{j=1}^d \sigma
_{i,j}^2\Biggr)  \int_0^t \mathrm{d}r \int_{\mathbb{R}^k} \mu(\mathrm{d}\xi)
\vert\mathcal{F}G(t-r)(\xi)\vert^2
\]
and it is well known (see \cite{t}) that
\begin{equation}
\label{fg}
\mathcal{F}G(t)(\xi)= \frac{\sin(t\Vert\xi\Vert
)}{\Vert\xi
\Vert}.
\end{equation}

Following \cite{d99} and \cite{pz}, we note that when $\mu$ is not
the null measure, the solution $u(t,x)$ of \eqref{1.1} is a random
vector, and the right-hand side of \eqref{1.3} is well defined if and
only if the following hypothesis is satisfied:
\renewcommand{\theequation}{H}
\begin{equation}\label{H}
0<\int_{\mathbb{R}^k
}\frac
{\mu(\mathrm{d}\xi)}{1+\Vert\xi\Vert^2} < \infty.
\end{equation}
In this case, the process $u$ given by \eqref{1.3} is a
natural example of an anisotropic Gaussian process, as considered in
\cite{x}.
Note that for the covariance density in (\ref{e1}), condition (\ref{H}) is satisfied when $\beta\in  \,]0,2\wedge k[$.

In Section \ref{s3} of this paper, we develop several results on
hitting probabilities that are related to those of \cite{dkn07}, but
are appropriate for studying the wave equation in \textit{all} spatial
dimensions. Indeed, the results of \cite{dkn07} were tailored to the
particularities of the heat equation in spatial dimension one, while
our results highlight the role of the spatial dimension and are
applicable to the stochastic wave equation. Theorem \ref{t3.1} gives a
lower bound on hitting probabilities;
Proposition \ref{p3.1} and Theorem \ref{t3.2} give upper bounds.
These three results apply to arbitrary stochastic processes, while
Theorem \ref{t3.3} gives a refinement of the upper bound in the case
of Gaussian processes. These results are used in Section \ref{s4}, but
will also be useful for studying nonlinear forms of \eqref{1.1}, which
is the subject of work currently in progress.

In Section \ref{s2}, we give simple conditions on a Gaussian process
$(X(t))$ that ensure an upper bound on the density function of
$(X(t),  X(s))$. This is related to a result in \cite{dn04}. The
upper bound is expressed in terms of the canonical metric of the
Gaussian process.

In Section \ref{s4}, the main effort is to obtain upper and lower
bounds on the behavior of the canonical metric associated with the
process $u$ (Proposition \ref{p4.2}). This is somewhat intricate for
the lower bounds, mainly because the expression for $E((u(t,x) -
u(s,y))^2)$ involves integrals of trigonometric functions and these are
not so easy to bound from below by positive quantities. Section \ref
{s4} ends by applying the results of Sections \ref{s3} and \ref{s2}
on hitting probabilities to obtain Theorems \ref{t4.1} and \ref
{t4.2}. These yield the following types of bounds:
\renewcommand{\theequation}{\arabic{equation}}
\setcounter{equation}{7}
\begin{equation}\label{e3}
c  \operatorname{Cap}_{d-{2(k+1)}/{(2-\beta)}}(A) \leq P\{u([t_0,T]\times
[-M,M]^k) \cap A\neq\varnothing\} \leq C  {\mathcal{H}}_{d-{2(k+1)}/{(2-\beta)}}(A),
\end{equation}
where Cap$_\gamma$ and $\mathcal{H}_\gamma$ denote capacity and
Hausdorff measure, respectively (their definitions are recalled in
Section \ref{s3}). We note that the same dimensions appear on both the
left- and right-hand sides of \eqref{e3}. This conclusion could also
have been deduced from Theorem 7.6 in \cite{x} or Theorem 2.1 in \cite
{BCX}, which contain general results on hitting probabilities for
anisotropic Gaussian processes. This is because our estimates on the
canonical metric of $u$ mentioned above, together with our Lemma \ref
{l2.1}, verify conditions (C1) and (C2) in these two references. We
also note that these estimates hint at the fact that condition (C3$'$) of
\cite{x} should be satisfied by $u$.

We recall that a point $z \in\mathbb{R}^d$ is \textit{polar} for $u$ if, for
all $t_0 >0$ and $M >0$,
\[
P\{z \in u([t_0,T]\times[-M,M]^k)\} = 0.
\]
Notice, as a consequence of \eqref{e3}, that if $d < 2(k+1)/(2-\beta
)$, then points \textit{are not} polar for $u$, while if $d >
2(k+1)/(2-\beta)$, then points \textit{are} polar for $u$. In the case
where $\beta$ is rational and $2(k+1)/(2-\beta)=d$ is an integer,
polarity of points in the critical dimension $d$ is an open problem.

As mentioned above, in a forthcoming work, we plan to extend these
results to systems of nonlinear stochastic wave equations with additive
noise, but without using Girsanov's theorem. It is a separate endeavor
to develop, using Malliavin calculus, the estimates needed for
multiplicative noise, as was done in \cite{dkn09} for the heat
equation. This will also make use of the results in Section~\ref{s3}.

\section{General results on hitting probabilities}
\label{s3}

Throughout this section, $V=\{v(x), x\in\mathbb{R}^m\}$, $m\in
\mathbb
{N}^*$, denotes an $\mathbb{R}^d$-valued stochastic process with continuous
sample paths.
We will fix a compact set $I\subset\mathbb{R}^m$ of positive Lebesgue
measure and consider an arbitrary Borel set $A\subset\mathbb{R}^d$.
Our aim is to give sufficient conditions on the stochastic process $V$
which lead to lower and upper bounds on the hitting probabilities
\[
P\{v(I)\cap A \ne\varnothing\}
\]
in terms of the capacity and the Hausdorff measure of $A$,
respectively, of a certain dimension. Here, $v(I)$ denotes the image of
$I$ under the (random)
map $x\mapsto v(x)$.

We now introduce some notation and recall the definition of capacity
and Hausdorff measure. For any $\gamma\in\mathbb{R}$, we define the
\textit{Bessel--Riesz kernels} by
\begin{equation}
\label{1.41}
K_\gamma(r)=
\cases{
r^{-\gamma} , &\quad  if $\gamma>0 $,\vspace*{2pt}\cr
\log\biggl(\dfrac{c}{r}\biggr) , & \quad if  $\gamma=0 $,\vspace*{2pt}\cr
1 , & \quad if  $\gamma<0 $,}
\end{equation}
where $c$ is a constant whose value will be specified later in the
proof of Lemma \ref{l3.1}.
Then, for every Borel set $A\subset\mathbb{R}^d$, we define $\mathcal
{P}(A)$ to be the set of probability measures on $A$. For \mbox{$\mu\in
\mathcal{P}(A)$}, we set
\[
\mathcal{E}_\gamma(\mu)=\int_{A} \int_{A}K_{\gamma}(\Vert
x-y\Vert) \mu(\mathrm{d}x)\mu(\mathrm{d}y).
\]
The \textit{Bessel--Riesz capacity} of a Borel set $A\subset\mathbb
{R}^d$ is
defined as follows:
\begin{equation}
\label{1.42}
\operatorname{Cap}_\gamma(A)=\Bigl[\inf_{\mu\in\mathcal{P}(A)}\mathcal
{E}_\gamma(\mu)\Bigr]^{-1}
\end{equation}
with the convention that $1/\infty=0$.

The $\gamma$-dimensional \textit{Hausdorff measure} of a Borel set
$A\subset\mathbb{R}^d$ is defined by $\mathcal{H}_\gamma(A)=\infty
$ if
$\gamma<0$, and for
\mbox{$\gamma\ge0$},
\begin{equation}
\label{1.43}
\mathcal{H}_\gamma(A)= \liminf_{\varepsilon\to0^+}\Biggl\{\sum
_{i=1}^\infty(2r_i)^\gamma\dvt A\subset\bigcup_{i=1}^\infty B_{r_i}(x_i),
\sup_{i\ge1}r_i\le\varepsilon\Biggr\}.
\end{equation}
Here, and throughout the paper, $B_r(x)$ denotes the open Euclidean
ball centered at $x$ and with radius $r$. Positive constants will be
denoted most often by
$C$ or $c$, although their value may change from one line to the next.
For a given subset $S\subset\mathbb{R}^n$ and $\nu>0$, we denote by
$S^{(\nu)}$ the $\nu$-enlargement of $S$.

We begin by studying the lower bound for $P\{v(I)\cap A \ne\varnothing\}$.

\begin{teorem}
\label{t3.1}
Fix $N>0$ and assume that the stochastic process $V$ satisfies the
following two hypotheses:
\begin{enumerate}[(2)]
\item[(1)] For any $x,y\in I$ with $x\neq y$, the vector $(v(x),
v(y))$ has density $p_{x,y}$ and there exist $\gamma, \alpha\in  \,]0,\infty[$
such that
\[
p_{x,y}(z_1,z_2) \le C\frac{1}{\Vert x-y\Vert^\gamma} \exp
\biggl(-\frac{c \Vert z_1-z_2\Vert^2}{\Vert x-y\Vert^\alpha}\biggr)
\]
for any $z_1,z_2\in[-N,N]^d$, where $C$ and $c$ are positive constants
independent of $x$ and $y$.
\item[(2)] One of the following two conditions holds:
\begin{enumerate}[(P$'$)]
\item[(P)] the density $p_x$ of $v(x)$ is continuous and bounded, and
$p_x(w)>0$ for any $x\in I$ and $w\in[-(N+1),N+1]^d$;
\item[(P$'$)] for any compact set $K\subset\mathbb{R}^d$ and any $x\in I$,
$\inf_{w\in K} p_x(w)\ge c_0>0$.
\end{enumerate}
\end{enumerate}
There then exists a positive and finite constant $c=c(N,\alpha,\gamma
,I,m)$ such that for all Borel sets $A\subset[-N,N]^d$,
\begin{equation}
\label{3.9}
P\{v(I)\cap A\neq\varnothing\} \ge c  \operatorname{Cap}_{({2}/{\alpha
})(\gamma-m)}(A).
\end{equation}
\end{teorem}

\begin{pf}
Without loss of generality, we may assume that $\operatorname{Cap}_{({2}/{\alpha})(\gamma-m)}(A)>0$,
otherwise there is nothing to prove. Under this assumption, we
necessarily have $\frac{2}{\alpha}(\gamma-m)<d$ and $A\ne\varnothing$
(see~\cite{k}, Appendix~C, Corollary 2.3.1, page 525).

First, assume that $A$ is a compact set. Following the scheme of the
proof of Theorem 2.1 in \cite{dkn07},
we consider three different cases.

\begin{longlist}
\item[\textit{Case 1}:] $\gamma-m<0$.
Let $z\in A$, $\varepsilon\in  \,]0,1[$ and set
\[
J_\varepsilon(z)=\frac{1}{(2\varepsilon)^d} \int_I \mathrm{d}x \, \mathbf{1}_{B_{\varepsilon}(z)}(v(x)).
\]
We will prove that $E(J_\varepsilon(z))\ge c_1$ and $E[(J_\varepsilon
(z))^2]\le c_2$ for some positive constants
$c_1$, $c_2$. With this, by using the Paley--Zygmund inequality (\cite
{k}, Chapter 3, Lemma 1.4.1) and noticing that
$\operatorname{Cap}_\beta(A)=1$ for $\beta<0$, we obtain
\begin{eqnarray*}
P\{J_\varepsilon(z) >0\}&\ge&\frac{[E (J_\varepsilon
(z))]^2}{E[(J_\varepsilon(z))^2]}
\ge C\\
&= &C \operatorname{Cap}_{({2}/{\alpha})(\gamma-m)}(A).
\end{eqnarray*}
However, $P\{J_\varepsilon(z) >0\}$ is bounded above by $P\{v(I)\cap
A^{(\varepsilon)}\ne\varnothing\}$. 
Since $A$ is compact and the trajectories of $v$ are continuous, by
letting $\varepsilon$ tend to $0$, we obtain (\ref{3.9}).

The lower bound for $E(J_\varepsilon(z))$ is a direct consequence of
assumption (2). To obtain the upper bound for
$E[(J_\varepsilon(z))^2]$, we first use the hypothesis (1) to obtain
\[
E[(J_\varepsilon(z))^2] \le C \int_I \mathrm{d}x \int_I \mathrm{d}y\, \frac{1}{\Vert
x-y\Vert^\gamma}.
\]
Let $\rho_0>0$ be such that $I\subset B_{{\rho_0}/{2}}(0)$. Fix
$x\in I$; after the change of variables $y\to x-y$
and by considering polar coordinates, we easily get
\[
E[(J_{\varepsilon}(z))^2] \le C \int_0^{\rho_0} \rho^{m-1-\gamma}
\,\mathrm{d}\rho.
\]
The last integral is bounded by a finite positive constant $c(m,\gamma
,I)$. Therefore, we obtain
$E[(J_\varepsilon(z))^2]\le c_2$.

\item[\textit{Case 2}:] $0<\frac{2}{\alpha}(\gamma-m)<d$.
Let $\mu\in\mathcal{P}(A)$. Let $g_\varepsilon=\frac
{1}{(2\varepsilon)^d}\mathbf{1}_{B_{\varepsilon}(0)}$ and
\[
J_{\varepsilon}(\mu) = \frac{1}{(2\varepsilon)^d}\int_I \mathrm{d}x \int
_{A} \mu(\mathrm{d}z) \mathbf{1}_{B_{\varepsilon}(0)}\bigl(v(x)-z\bigr)
= \int_I \mathrm{d}x \,(g_\varepsilon\ast\mu)(v(x)).
\]
Clearly, assumption (2) implies that $E(J_{\varepsilon}(\mu
))\ge c_1$ for a constant $c_1$ which does not depend on $\mu$
or $\varepsilon$. Moreover,
\[
E[(J_\varepsilon(\mu))^2] = \int_I \mathrm{d}x \int_I \mathrm{d}y \int_{\mathbb{R}^d}
\mathrm{d}z_1 \int_{\mathbb{R}^d} \mathrm{d}z_2\,(g_\varepsilon\ast\mu)(z_1)
(g_\varepsilon\ast\mu)(z_2) p_{x,y}(z_1,z_2).
\]
By hypothesis (1), Lemma \ref{l3.1} below and Theorem B.1 in \cite
{dkn07}, this is bounded by
\begin{eqnarray*}
&&C \int_{\mathbb{R}^d} \mathrm{d}z_1 \int_{\mathbb{R}^d}
\mathrm{d}z_2\,
(g_\varepsilon\ast
\mu)(z_1)(g_\varepsilon\ast\mu)(z_2) K_{({2}/{\alpha})(\gamma-m)}
(\Vert z_1-z_2\Vert)\\
&&\quad  =C \mathcal{E}_{({2}/{\alpha})(\gamma-m)}(g_\varepsilon
\ast\mu)\\
&&\quad \le C \mathcal{E}_{({2}/{\alpha})(\gamma-m)}(\mu).
\end{eqnarray*}
By choosing $\mu$ such that $\mathcal{E}_{({2}/{\alpha})(\gamma
-m)}(\mu)\le2/\operatorname{Cap}_{({2}/{\alpha})(\gamma-m)}(A)$, we obtain
\[
E[(J_\varepsilon(\mu))^2] \le\frac{C}{\operatorname{Cap}_{({2}/{\alpha})(\gamma-m)}(A)} ,
\]
and this yields (\ref{3.9}) by an argument similar to that used in
Case 1.

\item[\textit{Case 3}:] $\gamma-m=0$.
The proof is carried out in exactly the same way as for Case 2, by
applying Theorem B.2 in \cite{dkn07} instead of Theorem B.1.

Now, let $A$ be a Borel set included in $[-N,N]^d$. It is well known that
\begin{equation}
\label{compact}
\operatorname{Cap}_{\beta}(A) = \sup_{F\subset A,  F  \mathrm{\ compact}}\operatorname{Cap}_{\beta}(F)
\end{equation}
(see, for instance, Chapter 3 of \cite{dm}). Therefore, for any
compact set $F\subset A$, we have
\[
P\{v(I)\cap A \ne\varnothing\}\ge P\{v(I)\cap F \ne\varnothing\}\ge c \operatorname{Cap}_{({2}/{\alpha})(\gamma-m)}(F).
\]
This yields (\ref{3.9}) by taking the supremum over such $F$ and using
(\ref{compact}).
\end{longlist}

The proof of the theorem is thus complete.
\end{pf}

In order to end the study of lower bounds, we prove a technical lemma
which was used in the proof of Theorem \ref{t3.1} to relate joint
densities to Bessel--Riesz kernels.

\begin{lem}
\label{l3.1}
Fix $\alpha, \gamma\in  \,]0,\infty[$.
There exists a constant
$C:=C(N,\alpha,\gamma,I,m)$ such that for any $a\in  \,]{-}N,N[$,
\begin{equation}
\label{3.7.1}
\int_I \mathrm{d}x \int_I \mathrm{d}y\, \frac{1}{\Vert x-y\Vert^\gamma}\exp
\biggl(-\frac{a^2}{\Vert x-y\Vert^\alpha}\biggr) \le C K_{({2}/{\alpha})(\gamma-m)}(a).
\end{equation}
\end{lem}

\begin{pf}
Fix $\rho_0>0$ such that $I\subset B_{{\rho
_0}/{2}}(0)$. Fix $x\in I$ and consider the change of variables
$z=a^{-{2}/{\alpha}}(x-y)$. Denoting by ${\mathcal{I}}$ the
left-hand side of (\ref{3.7.1}), we have
\[
{\mathcal{I}}\le C(I)  a^{-({2}/{\alpha})(\gamma-m)} \int
_{B_{{\rho_0}/{a^{2/\alpha}}}(0)} \mathrm{d}z \frac{1}{\Vert z\Vert
^\gamma} \exp\biggl(-\frac{1}{\Vert z\Vert^\alpha}\biggr).
\]

Let
\[
{\mathcal{J}}= \int_{B_{{\rho_0}/{a^{2/\alpha}}}(0)} \mathrm{d}z \,\frac
{1}{\Vert z\Vert^\gamma} \exp\biggl(-\frac{1}{\Vert z\Vert^\alpha
}\biggr).
\]
Using polar coordinates, we have ${\mathcal{J}}={\mathcal
{J}}_1+{\mathcal{J}}_2$, where
\begin{eqnarray*}
{\mathcal{J}}_1&=&\int_0^{{\rho_0}/{N^{2/\alpha}}} \mathrm{d}\rho\, \rho
^{m-1-\gamma} \exp\biggl(-\frac{1}{\rho^\alpha}\biggr),\\
{\mathcal{J}}_2&=&\int_{{\rho_0}/{N^{2/\alpha}}}^{{\rho
_0}/{a^{2/\alpha}}} \mathrm{d}\rho \,\rho^{m-1-\gamma} \exp\biggl(-\frac
{1}{\rho^\alpha}\biggr).
\end{eqnarray*}
Clearly, ${\mathcal{J}}_1\le C(\rho_0,N)$. In order to study
${\mathcal{J}}_2$, we bound the exponential by $1$ and consider three
different cases.

\begin{longlist}
\item[\textit{Case 1}.] If $m-\gamma<0$, then
\[
{\mathcal{J}}_2\le(\gamma-m)^{-1}\biggl(\frac{\rho
_0}{N^{{2}/{\alpha}}}\biggr)^{m-\gamma}\le C(N,\alpha,\gamma,
\rho_0,m).
\]
\item[\textit{Case 2}.] If $m-\gamma>0$, then
\[
{\mathcal{J}}_2\le(m-\gamma)^{-1}\biggl(\frac{\rho
_0}{a^{{2}/{\alpha}}}\biggr)^{m-\gamma}\le C(\gamma, \rho
_0,m)a^{({2}/{\alpha})(\gamma-m)}.
\]
\item[\textit{Case 3}.] If $m-\gamma=0$, then
\[
{\mathcal{J}}_2 \le\frac{2}{\alpha}\log\biggl(\frac{N}{a}\biggr).
\]
\end{longlist}
Since ${\mathcal{I}}\le C(I) a^{-({2}/{\alpha})(\gamma
-m)}{\mathcal{J}}$, we reach the conclusion using the definition of
$K_\beta(a)$ for
\mbox{$\beta,a\in\mathbb{R}$} (see (\ref{1.41}); in the case where $m -
\gamma=
0$, the constant $c$ in \eqref{1.41} must be chosen sufficiently large).
\end{pf}

We now study upper bounds for the hitting probabilities. For this,
we fix $\delta>0$, $\varepsilon\in  \,]0,1[$, $j_1, \ldots, j_m\in
\mathbb{Z}$, and set $j=(j_1, \ldots, j_m)$ and
\begin{equation}
\label{r}
R_{j}^\varepsilon= \prod_{l=1}^m[j_l \varepsilon^{{1}/{\delta}},( j_l+1) \varepsilon^{{1}/{\delta}}].
\end{equation}
%

The next statement is an extension to higher dimensions of Theorem 3.1
in \cite{dkn07}.

\begin{prop}
\label{p3.1}
Let $D\subset\mathbb{R}^d$ and $\gamma>0$. 
We assume that
there exists a positive constant $c$ such that, for all small
$\varepsilon\in  \,]0,1[$, $z\in D^{(1)}$
and any set $R_{j}^\varepsilon$ such that $R_{j}^\varepsilon\cap I\ne
\varnothing$,
\begin{equation}
\label{3.10}
P\{v(R_{j}^\varepsilon)\cap B_\varepsilon(z) \neq\varnothing
\}\le c  \varepsilon^\gamma.
\end{equation}
There then exists a positive constant $C$ such that for any Borel set
$A\subset D$,
\begin{equation}
\label{3.11}
P\{v(I)\cap A \neq\varnothing\}\le C \mathcal{H}_{\gamma
-{m}/{\delta}}(A).
\end{equation}
%
\end{prop}

\begin{pf}
We suppose that $\gamma-\frac{m}{\delta}\ge0$,
otherwise $\mathcal{H}_{\gamma-{m}/{\delta}}(A)=\infty$ and
therefore (\ref{3.11}) obviously holds.
Clearly, by the additive property of probability,
\[
P\{v(I)\cap B_\varepsilon(z) \neq\varnothing\}\le\sum
_{j: R_{j}^\varepsilon\cap I\neq\varnothing}P\{
v(R_{j}^\varepsilon)\cap B_\varepsilon(z) \neq\varnothing\}
\]
for any $\varepsilon>0$.
Since $I$ is bounded, the number of terms in the sum on the right-hand
side of this inequality is bounded by a multiple of $\varepsilon
^{-{m}/{\delta}}$. Hence,
\[
P\{v(I)\cap B_\varepsilon(z) \neq\varnothing\}\le C
\varepsilon^{-{m}/{\delta}}P\{v(R_{j}^\varepsilon)\cap
B_\varepsilon(z) \neq\varnothing\}.
\]
Using (\ref{3.10}), we then obtain
\begin{equation}
\label{3.12}
P\{v(I)\cap B_\varepsilon(z) \neq\varnothing\}\le C
\varepsilon^{\gamma-{m}/{\delta}}.
\end{equation}
This yields (\ref{3.11}) by a \textit{covering argument}, as shown in the
proof of Theorem 3.1 in \cite{dkn07}.
For the sake of completeness, we sketch this argument.

Fix $\varepsilon\in  \,]0,1[$ sufficiently small and consider a
sequence of open balls $(B_n, n\ge1)$ with respective radii $r_n\in  \,]0,\varepsilon]$,
such that $B_n\cap A\ne\varnothing$, $A\subset\bigcup_{n\ge1}B_n$ and
\[
\sum_{n\ge1}(2r_n)^{\gamma-{m}/{\delta}}\le\mathcal
{H}_{\gamma-{m}/{\delta}}(A)+\varepsilon.
\]
Then, by (\ref{3.12}),
\begin{eqnarray*}
P\{v(I)\cap A \neq\varnothing\}&\le&\sum_{n\ge1}P\{
v(I)\cap B_n\neq\varnothing\}\\
&\le& C\sum_{n\ge1} (2r_n)^{\gamma-{m}/{\delta}}\\
&\le& C\bigl(\mathcal{H}_{\gamma-{m}/{\delta}}(A)+\varepsilon
\bigr).
\end{eqnarray*}
Finally, we let $\varepsilon\downarrow0$ to complete the proof.
\end{pf}

In the next theorem, we give sufficient conditions on the process $V$
for the assumptions of Proposition \ref{p3.1} to be satisfied and
therefore to
ensure (\ref{3.11}).

\begin{teorem}
\label{t3.2}
Let $D\subset\mathbb{R}^d$. Assume that the stochastic process $V$ satisfies
the following two conditions:
%
\begin{enumerate}[(2)]
\item[(1)] for any $x\in\mathbb{R}^m$, the random vector $v(x)$ has density
$p_x$ and
\[
\sup_{z\in D^{(2)}} \sup_{x\in I ^{(1)}} p_{x}(z) \le C ;
\]
\item[(2)] there exists $\delta\in  \,]0,1]$ and a constant $C$ such
that for any $q\in[1,\infty[$, $x,y\in I^{(1)}$,
%
\[
E\bigl(\Vert v(x)-v(y)\Vert^q\bigr) \le C\Vert
x-y\Vert^{q\delta}.
\]
\end{enumerate}
%
For any $\gamma\in  \,]0,d[$, inequality \textup{(\ref{3.10})} then holds
and, consequently, for every Borel set $A\subset D$,
\begin{equation}
\label{3.13}
P\{v(I)\cap A \neq\varnothing\}\le C\mathcal{H}_{\gamma
-{m}/{\delta}}(A).
\end{equation}
\end{teorem}

\begin{pf}
We keep the notation of Proposition \ref{p3.1} and write
$x_j^\varepsilon=(j_l \varepsilon^{{1}/{\delta}}, l=1,\ldots
,m)$. For any $z\in D^{(1)}$ and $R_j^\varepsilon$ such that
$R_j^\varepsilon\cap I \ne\varnothing$, set
\[
Y_j^\varepsilon= \Vert v(x_j^\varepsilon)-z\Vert,\qquad  Z_j^\varepsilon= \sup_{x\in R_j^\varepsilon}\Vert v(x) -
v(x_j^\varepsilon)\Vert.
\]
By applying the version of Kolmogorov's criterion as it is stated in
\cite{ry}, Theorem 2.1, page 26, using assumption (2), we obtain
\[
E((Z_{j}^\varepsilon)^q)\le C\Vert
x-x_j^\varepsilon\Vert^{\alpha q}
\]
for any $q\in[1,\infty[$ and $\alpha\in  \,]0,\delta-\frac{m}{q}[$.
Hence,
\begin{equation}
\label{3.14}
E(( Z_{j}^\varepsilon)^q)\le C \varepsilon
^{\gamma_0q}
\end{equation}
with $\gamma_0<1-\frac{m}{q\delta}$.\vspace*{1.5pt}

Let $\gamma\in  \,]0,d[$. We first prove that
\begin{equation}
\label{3.15}
P\bigl\{Z_{j}^\varepsilon\ge\tfrac{1}{2}Y_{j}^\varepsilon\bigr\} \le C
\varepsilon^\gamma.
\end{equation}
For this, we consider the decomposition
\[
P\bigl\{Z_{j}^\varepsilon\ge\tfrac{1}{2}Y_{j}^\varepsilon\bigr\} \le P\{
Y_{j}^\varepsilon\le\varepsilon^{{\gamma}/{d}}\} + P\bigl\{
Z_{j}^\varepsilon\ge\tfrac{1}{2}\varepsilon^{{\gamma}/{d}}\bigr\}
\]
and then give upper bounds for each term on the right-hand side.

Clearly, from the boundedness of the density stated in assumption (1),
\[
P\{Y_{j}^\varepsilon\le\varepsilon^{{\gamma}/{d}}\} \le C
\varepsilon^\gamma
\]
and by Markov's inequality, along with (\ref{3.14}),
we have
\[
P\bigl\{Z_{j}^\varepsilon\ge\tfrac{1}{2}\varepsilon^{{\gamma}/{d}}\bigr\}
\le C \varepsilon^{q(\gamma_0-{\gamma}/{d})}.
\]
Therefore,
\[
P\bigl\{Z_{j}^\varepsilon\ge\tfrac{1}{2}Y_{j}^\varepsilon\bigr\} \le C
\bigl(\varepsilon^\gamma+ \varepsilon^{q(\gamma_0-{\gamma
}/{d})}\bigr)
\]
for any $\gamma_0<1-\frac{m}{q\delta}$.
Since $\gamma\in  \,]0,d[$, we can choose $\gamma_0 < 1$ and $q$
arbitrarily large such that $\frac{\gamma}{d}<\gamma_0<1-\frac
{m}{q\delta}$. Hence, we obtain (\ref{3.15}).

If $v(R_{j}^\varepsilon)\cap B_\varepsilon(z)\neq\varnothing$, then
$Y_{j}^\varepsilon\le\varepsilon+ Z_{j}^\varepsilon$. Therefore,
\begin{eqnarray*}
P\{v(R_{j}^\varepsilon)\cap B_\varepsilon(z)\neq\varnothing\}&\le& P\{
Y_{j}^\varepsilon\le\varepsilon+ Z_{j}^\varepsilon\}\\
&\le& P\bigl\{Z_{j}^\varepsilon\ge\tfrac{1}{2}Y_{j}^\varepsilon\bigr\} + P\{
Y_{j}^\varepsilon\le2\varepsilon\}\\
&\le& C(\varepsilon^\gamma+\varepsilon^d)\\
& \le& C \varepsilon^\gamma
\end{eqnarray*}
since $\gamma\in  \,]0,d[$.
This proves (\ref{3.10}) for any $\gamma\in  \,]0,d[$. By Proposition
\ref{p3.1}, we obtain (\ref{3.13}).
\end{pf}

The remainder of this section is devoted to extending the validity of
(\ref{3.13}) to $\gamma=d$
in the case where $V$ belongs to a particular class of Gaussian
processes. For this class, we will prove that,
instead of (\ref{3.14}), the following, stronger, property holds:

For any $\varepsilon\in  \,]0,1[$, each $j\in\mathbb{Z}^m$ with
$R_j^\varepsilon\cap I\ne\varnothing$ and every $q\in[1,\infty[$,
there exists $C>0$ such that
\begin{equation}
\label{moments}
E\Bigl(\sup_{x\in R_j^\varepsilon}\Vert v(x)-v(x_j^\varepsilon
)\Vert^q\Bigr) \le C \varepsilon^q.
\end{equation}
We will then show that $P\{v(R_j^\varepsilon)\cap
B_\varepsilon(z)\ne\varnothing\}\le C \varepsilon^d$ (see Theorem
\ref{t3.3} below). Together with Proposition \ref{p3.1}, this will
yield the desired improvement.

We first give a sufficient condition which applies to
arbitrary continuous stochastic pro\-cesses~$V$.

\begin{lem}
\label{l3.2}
Let $\nu\in  \,]0,1]$. Suppose that for any $\varepsilon\in  \,]0,1[$
sufficiently small,
\begin{equation}
\label{3.16}
E\biggl(\int_{B_\varepsilon(x)} \mathrm{d}y \int_{B_\varepsilon(x)} \mathrm{d}\bar
{y}\, \biggl[\exp\biggl\{\frac{\Vert v(y)-v(\bar y)\Vert}{\Vert y-\bar
y\Vert^\nu}\biggr\}\biggr]\biggr) \le C \varepsilon^{2m} ,
\end{equation}
where $C$ is a positive constant.
Let $S^\nu_\varepsilon(x)=\{y\in\mathbb{R}^m\dvt \Vert x-y\Vert\le
\varepsilon^{{1}/{\nu}}\}$.

Then, for any $q\in[1,\infty[$, there exists $\bar C>0$ such that for all
small $\varepsilon>0$,
\begin{equation}
\label{3.17}
E\Bigl(\sup_{y\in S^\nu_\varepsilon(x)} \Vert v(x)-v(y)\Vert
^q\Bigr) \le\bar C \varepsilon^q.
\end{equation}
\end{lem}

\begin{pf}
By (\ref{3.16}), $B(\omega)<\infty$ a.s., where
\[
B(\omega)= \int_{S^\nu_\varepsilon(x)} \mathrm{d}y \int_{S^\nu_\varepsilon
(x)} \mathrm{d}\bar{y}\,\biggl [\exp\biggl\{\frac{\Vert v(y)(\omega)-v(\bar
y)(\omega)\Vert}{\Vert y-\bar y\Vert^\nu}\biggr\}\biggr].
\]
We apply the Garsia--Rodemich--Rumsey lemma (see \cite{sv}, exercise
2.4.1, page 60) to the functions
$\psi(x)= \mathrm{e}^x-1$, $p(x)=x^\nu$ and functions $f\dvtx S^\nu_\varepsilon
(x)\subset\mathbb{R}^m\to\mathbb{R}^d$ given by the sample paths of
the process $V$
restricted to the parameter set $S^\nu_\varepsilon(x)$, to obtain
\[
\Vert v(x)-v(y)\Vert\le8 \int_0^{2\Vert x-y\Vert} \psi^{-1}
\biggl(\frac{C_1 B(\omega)}{u^{2m}}\biggr) \nu u^{\nu-1} \,\mathrm{d}u ,
\]
where $C_1$ is a positive constant which depends only on $m$.
Consequently, for any $q\in[1,\infty[$,
\[
E\Bigl(\sup_{y\in S^\nu_\varepsilon(x)}\Vert v(x)-v(y)\Vert
^q\Bigr)\le8 E\biggl(\int_0^{2 \varepsilon^{{1}/{\nu}}}\psi
^{-1}\biggl(\frac{C_1 B(\omega)}{u^{2m}}\biggr) \nu u^{\nu-1}\, \mathrm{d}u
\biggr)^q.
\]
We note that since $\psi^{-1}(x)= \ln(1+x)$ is an increasing function
on $[0,\infty)$, the constant $C_1$ above can be taken arbitrarily
large. In the sequel, we will fix $q\in[1,\infty[$ and assume that
$C_1\ge(\mathrm{e}^{q-1}-1) C_2^{-1}2^{2m}$, where $C_2$ is the
square of the volume of the unit ball in $\mathbb{R}^m$. Then
\[
B(\omega)\ge C_2\varepsilon^{{2m}/{\nu}}\ge\frac{\mathrm{e}^{q-1}-1}{C_1}u^{2m}
\]
for any $u\in[0,2\varepsilon^{{1}/{\nu}}]$.

Jensen's inequality applied first to the convex function $\varphi
_1(x)= x^q$, $x\in\mathbb{R}$, and the integral with respect to the measure
$\mu(\mathrm{d}u)=u^{\nu-1}\, \mathrm{d}u$, and then to the concave
function $\varphi_2(x)= \ln^q(1+x)$, $x\in[\mathrm{e}^{q-1}-1,\infty[$, and
to the expectation operator, yields
\begin{eqnarray*}
E\Bigl(\sup_{y\in S^\nu_\varepsilon(x)}\Vert v(x)-v(y)\Vert
^q\Bigr)&\le&8 \varepsilon^{q-1}\int_0^{2\varepsilon^{{1}/{\nu}}} E\biggl[\ln
^q\biggl(1+\frac{C_1B(\omega)}{u^{2m}}\biggr)\biggr] u^{\nu-1} \,\mathrm{d}u\\
&\le& C\varepsilon^{q-1}\int_0^{2 \varepsilon^{{1}/{\nu}}} \ln
^q\biggl(1+\frac{C_3\varepsilon^{{2m}/{\nu}}}{u^{2m}}\biggr)
\nu u^{\nu-1} \,\mathrm{d}u
\end{eqnarray*}
with $C_3=C_1 C$.
With the change of variable $u\to\frac{u^\nu}{\varepsilon}$, we have
\begin{eqnarray*}
\int_0^{2\varepsilon^{{1}/{\nu}}} \ln^q\biggl(1+\frac{C_3
\varepsilon^{{2\mu}/{\nu}}}{u^{2m}}\biggr) \nu u^{\nu-1}\, \mathrm{d}u
&=&\varepsilon\int_0^{2^\nu} \ln^q\biggl(1+\frac{C_3}{w^{
{2m}/{\nu}}}\biggr) \,\mathrm{d}w\\
&=& \bar C \varepsilon.
\end{eqnarray*}
This proves (\ref{3.17}).
\end{pf}

We can now sharpen the result of Theorem \ref{t3.2} in the case of
Gaussian processes.

\begin{teorem}
\label{t3.3}
Assume that the stochastic process $V=\{v(x),  x\in\mathbb{R}^m\}$ is
continuous, Gaussian and centered, with independent, identically
distributed components
$\{v_i(x),  x\in\mathbb{R}^m\}$, $i=1,\ldots,d$,
and $\inf_{x\in I^{(1)}}{\operatorname{Var}} (v_1(x))>0$. Fix $\delta\in  \,]0,1]$ and suppose that for any $\varepsilon>0$ small enough and any
$R_j^\varepsilon$
(defined in \textup{(\ref{r})}) such that $R_j^\varepsilon\cap I \ne\varnothing$,
we have
\begin{equation}
\label{m}
E\biggl(\int_{R_j^\varepsilon} \mathrm{d}y \int_{R_j^\varepsilon}
\mathrm{d}\bar{y}\,
\biggl[\exp\biggl\{\frac{\Vert v(y)-v(\bar y)\Vert}{\Vert y-\bar
y\Vert^\delta}\biggr\}\biggr]\biggr) \le C \varepsilon^{
{2m}/{\delta}}.
\end{equation}
Then, for every $z\in\mathbb{R}^d$ and $R_j^\varepsilon$ as before,
\begin{equation}
\label{3.18}
P\{v(R_j^\varepsilon)\cap B_\varepsilon(z)\ne\varnothing\}\le C
\varepsilon^d.
\end{equation}
Consequently, for any Borel set $A\subset\mathbb{R}^d$,
\begin{equation}
\label{3.18.1}
P\{v(I)\cap A\ne\varnothing\}\le C \mathcal{H}_{d-{m}/{\delta
}}(A).
\end{equation}
\end{teorem}

\begin{pf}
By Lemma \ref{l3.2}, assumption (\ref{m}) implies (\ref
{moments}).
We use this property and adapt the proof of Proposition 4.4 of \cite
{dkn07}. First, for any $z\in\mathbb{R}^d$, we write
\[
P\{v(R_j^\varepsilon)\cap B_\varepsilon(z)\ne\varnothing\}= P\Bigl\{
\inf_{x\in R_j^\varepsilon}\Vert v(x)-z\Vert\le\varepsilon\Bigr\}.
\]
Next, we write the condition $\Vert v(x)-z\Vert\le\varepsilon$ in
terms of two independent random variables, as follows.
Set
\[
c_j^\varepsilon(x)= \frac{E(v_1(x)v_1(x_j^\varepsilon)
)}{{\operatorname{Var}}(v_1(x_j^\varepsilon))}.
\]
Because $V$ is a Gaussian process,
\[
E(v(x)\vert v(x_j^\varepsilon)) = c_j^\varepsilon(x)
v(x_j^\varepsilon).
\]
Set
\[
Y_j^\varepsilon=\inf_{x\in R_j^\varepsilon}\Vert c_j^\varepsilon
(x)v(x_j^\varepsilon)-z\Vert ,\qquad
Z_j^\varepsilon=\sup_{x\in R_j^\varepsilon}\Vert
v(x)-c_j^\varepsilon(x)v(x_j^\varepsilon)\Vert.
\]
Again, because $V$ is a Gaussian process, these two random variables
are independent and
\begin{equation}
\label{3.19}
P\Bigl\{\inf_{x\in R_j^\varepsilon}\Vert v(x)-z\Vert\le\varepsilon\Bigr\}\le
P\{Y_j^\varepsilon\le\varepsilon+Z_j^\varepsilon\}.
\end{equation}
Our next aim is to prove that for any $r\ge0$,
\begin{equation}
\label{3.20}
P(Y_j^\varepsilon\le r)\le C r^d .
\end{equation}
For this, we first note that by independence of the components of $V$,
\[
P(Y_j^\varepsilon\le r)\le\prod_{i=1}^d P(G_{j,i}^\varepsilon) ,
\]
where
\[
G_{j,i}^\varepsilon= \Bigl\{\inf_{x\in R_j^\varepsilon}|c_j^\varepsilon
(x)v_i(x_j^\varepsilon)-z_i|\le r\Bigr\}.
\]
By setting $e_j^\varepsilon= \inf_{x\in R_j^\varepsilon
}c_j^\varepsilon(x)$, we have
\[
P( G_{j,i}^\varepsilon)\le P\bigl(v_i(x_j^\varepsilon)\in B_{{r}/{e_j^\varepsilon}}(z)\bigr).
\]
Since $V$ is centered and $\inf_{x\in I^{(1)}}{\rm Var}(v_1(x))>0$ by
hypothesis, Schwarz's inequality and (\ref{moments}) yield
\begin{eqnarray}\label{3.21}
|c_j^\varepsilon(x)-1|&=&\frac{\vert E[v_1(x_j^\varepsilon
)(v_1(x)-v_1(x_j^\varepsilon))]\vert}{{\operatorname{Var}}(v_1(x_j^\varepsilon))}\nonumber\\
&\le& C \biggl(\frac{E([v_1(x)-v_1(x_j^\varepsilon)
]^2)}{{\operatorname{Var}}(v_1(x_j^\varepsilon))}\biggr)^{{1}/{2}}\\
&\le& C \varepsilon\nonumber
\end{eqnarray}
for any $x\in R_j^\varepsilon$.
This implies that $\frac{r}{e_j^\varepsilon}\le C r$ and since the
density of $v_i(x_j^\varepsilon)$ is bounded, we get
\[
P\bigl(v_i(x_j^\varepsilon)\in B_{{r}/{e_j^\varepsilon
}}(z)\bigr)\le C r.
\]
Therefore, (\ref{3.20}) holds.

By \eqref{3.20} and the independence of $Y_j^\varepsilon$ and
$Z_j^\varepsilon$,
\[
P\{Y_j^\varepsilon\le\varepsilon+Z_j^\varepsilon\}\le C E
[(\varepsilon+Z_j^\varepsilon)^d].
\]
Consider the decomposition $Z_j^\varepsilon=Z_j^{\varepsilon
,1}+Z_j^{\varepsilon,2}$, where
\[
Z_j^{\varepsilon,1}=\sup_{x\in R_j^\varepsilon}\Vert
v(x)-v(x_j^\varepsilon)\Vert ,\qquad
Z_j^{\varepsilon,2}=\Vert v(x_j^\varepsilon)\Vert\sup_{x\in
R_j^\varepsilon}|1-c_j^\varepsilon(x)|.
\]
By (\ref{moments}), we have $E(|Z_j^{\varepsilon,1}
|^d)\le C \varepsilon^d$.
Moreover, by (\ref{3.21}) and (\ref{m}),
\[
E(\Vert Z_j^{\varepsilon,2}\Vert^d)\le
C\varepsilon^d E(\Vert v(x_j^\varepsilon)\Vert^d)\le C
\varepsilon^d.
\]
This completes the proof of (\ref{3.18}). Finally, (\ref{3.18.1})
follows from Proposition \ref{p3.1}.
\end{pf}

\section{Joint densities of Gaussian processes}
\label{s2}
Consider a Gaussian family of centered, $\mathbb{R}^d$-valued random vectors,
indexed by a compact metric space $(\mathbb{T},d)$, that we denote by
$X=(X_t,  t\in\mathbb{T})$. We suppose that the component processes
$(X_t^i,  t\in\mathbb{T})$, $i=1,\ldots,d$, are independent.
We also assume mean square continuity, that is,
by letting
\[
\delta(s,t)= \bigl(E(\Vert X_t-X_s\Vert^2)
\bigr)^{{1}/{2}}
\]
denote the canonical (pseudo)-metric associated with $X$, we have
$\delta(s,t)\to0$ as $d (s,t)\to0$.

Let $p_{s,t}(z_1,z_2)$ denote the joint density of $(X_s,X_t)$ at
$(z_1,z_2)\in\mathbb{R}^{2d}$. The purpose of this section is to establish
upper bounds of exponential type for $p_{s,t}(z_1,z_2)$. We note that
these conditions and, in particular, condition (c) below, are easily
verified in many examples.

\begin{prop}
\label{p2.1}
Suppose that:
\begin{enumerate}[(c)]
\item[(a)] $\sigma_{t,i}^2 := \operatorname{Var}(X_t^i)>0$ for any
$i=1,\ldots,d$ and all $t\in\mathbb{T}$;
\item[(b)] $\operatorname{Corr}(X_s^i,X_t^i)<1$ for any $i=1,\ldots,d$ and
$s,t\in\mathbb{T}$ with $s\neq t$;
\item[(c)] there exists $\eta>0$ and a positive constant $C>0$ such
that for all $s,t\in\mathbb{T}$,
\[
\sup_{i\in\{1,\ldots,d\}}|\sigma_{t,i}^2-\sigma_{s,i}^2| \le
C(\delta(s,t))^{1+\eta}.
\]
\end{enumerate}
Fix $M>0$. There then exists $C>0$ such that for all $s,t\in\mathbb
{T}$ with $s \neq t$ and $z_1,z_2\in[-M,M]^d$,
\[
p_{s,t}(z_1,z_2) \le\frac{C}{(\delta(s,t))^d} \exp
\biggl(-\frac{c\Vert z_1-z_2\Vert^2}{(\delta(s,t))^2}\biggr)
\]
for some positive and finite constants $C$ and $c$.
\end{prop}

\begin{pf}
Note that (a), (b) and the independence of the components
yield the existence of $p_{s,t}$.

Fix $i=1,\ldots, d$ and denote by $p_{s,t}^i(z_1,z_2)$,
$p^i_{t|s}(\cdot|z_2)$ and $p^i_s(\cdot)$ the joint density of
$(X^i_s,X^i_t)$ at\vspace*{-1pt} $(z_1,z_2)$,
the conditional density of $X_t^i$ given $X_s^i=z_2$ and the marginal
density of $X^i_s$, respectively. It is well known (linear regression) that
\[
p^i_{t|s}(z_1|z_2)=\frac{1}{\tau_{s,t}\sqrt{2\curpi}} \exp
\biggl(-\frac{|z_1-m_{s,t}z_2|^2}{2\tau_{s,t}^2}\biggr),
\]
where
\[
\tau_{s,t}^2=\sigma_t^2(1-\rho_{s,t}^2) , \qquad  \rho_{s,t}=\frac
{\sigma_{s,t}}{\sigma_s\sigma_t} ,\qquad   m_{s,t}= \frac{\sigma
_{s,t}}{\sigma_s^2} ,\qquad   \sigma_{s,t}=E(X_s^iX_t^i)
\]
and, for the sake of simplicity, we have omitted the index $i$.
Since
\[
p_{s,t}^i(z_1,z_2)= p^i_{t|s}(z_1\vert z_2)p^i_s(z_2) ,
\]
the triangle inequality, along with the elementary bound $(a-b)^2\ge
\frac{1}{2}a^2-b^2$, yields
\begin{eqnarray*}
p_{s,t}^i(z_1,z_2)&\le&\frac{1}{2\curpi\sigma_s\tau_{s,t}}\exp
\biggl(-\frac{|z_1-z_2|^2}{4\tau_{s,t}^2}\biggr)\\
&&{} \times\exp\biggl(\frac{|z_2|^2|1-m_{s,t}|^2}{2\tau
_{s,t}^2}\biggr)
\exp\biggl(-\frac{|z_2|^2}{2\sigma_s^2}\biggr).
\end{eqnarray*}
By hypotheses (a) and (c), $s\mapsto\sigma_s^2$ is bounded above and
below by positive constants, therefore, for $z_2\in[-M,M]$,
\[
p_{s,t}^i(z_1,z_2)\le\frac{C}{\tau_{s,t}}\exp\biggl(-\frac
{|z_1-z_2|^2}{4\tau_{s,t}^2}\biggr) \exp\biggl(\frac
{M^2|1-m_{s,t}|^2}{2\tau_{s,t}^2}\biggr).
\]
The conclusion now follows from Lemma \ref{l2.1} below and the
independence of the components of~$X$.
\end{pf}

\begin{lem}
\label{l2.1}
With the same assumptions and notation as in Proposition \textup{\ref{p2.1}},
there exist constants $0<c_1<c_2<\infty$ such that for all $s,t\in
\mathbb{T}$,
\begin{enumerate}[(2)]
\item[(1)] $c_1 \delta(s,t)\le\tau_{s,t}\le c_2 \delta(s,t)$;
\item[(2)] $|1-m_{s,t}|\le c_2 \delta(s,t)$.
\end{enumerate}
\end{lem}

\begin{pf}
A simple calculation gives
\begin{equation}
\label{2.1}
\sigma_t^2\sigma_s^2-\sigma_{s,t}^2
=\tfrac{1}{4} [\delta(s,t)^2 - (\sigma_t-\sigma_s)^2
][(\sigma_t+\sigma_s)^2 - \delta(s,t)^2]
\end{equation}
(see \cite{mt}, equation (3.1)). Therefore, by hypothesis (c) of
Proposition \ref{p2.1},
\[
1-\rho_{s,t}^2 \le\frac{C}{\sigma_s^2\sigma_t^2} \delta(s,t)^2.
\]
From assumptions (a) and (c), it follows that there is a positive
constant $c_2<\infty$ such that for all $s,t\in\mathbb{T}$,
\[
\tau_{s,t}^2\le c_2^2  \delta(s,t)^2 ,
\]
which proves the upper bound in assertion (1).

For the lower bound in (1), we note that for $s$ near $t$, the second
factor on the right-hand side of (\ref{2.1}) is bounded below since
$\delta(s,t)\to0$
as $d(s,t)\to0$. Further, by hypotheses (a) and (c),
we have
\begin{eqnarray*}
\delta(s,t)^2 - (\sigma_t-\sigma_s)^2& =& \delta(s,t)^2 - \frac
{(\sigma_t^2-\sigma_s^2)^2}{(\sigma_t+\sigma
_s)^2}\\
&\ge&\delta(s,t)^2-\tilde c_1 \delta(s,t)^{2+2\eta}\\
&\ge& c_1\delta(s,t)^2
\end{eqnarray*}
for $s$ near $t$. This proves the lower bound in (1) when $\delta
(s,t)$ is sufficiently small.

In order to extend this inequality to all $s,t\in\mathbb{T}$, it
suffices to observe that by hypothesis (b),
\[
\sigma_t^2\sigma_s^2 - \sigma_{s,t}^2 > 0
\]
if $s\neq t$, and by hypothesis (c), for $\varepsilon> 0$, there
exists $c^\prime>0$ such that $\sigma_t^2\sigma_s^2-\sigma_{s,t}^2
>c^\prime$ for
$\delta(s,t)\ge\varepsilon$. This proves the lower bound in
assertion (1).

In order to prove assertion (2), observe that
\[
|1-m_{s,t}| = \frac{|\sigma_s^2 - \sigma_{t,s}|}{\sigma_s^2}
\]
and
\begin{eqnarray*}
|\sigma_s^2 - \sigma_{t,s}|& = &\big\vert\delta(s,t)^2 + E
\bigl((X_s-X_t)X_t\bigr)\big\vert\\
&\le&\delta(s,t)^2+ \delta(s,t) \sigma_t \le c  \delta(s,t).
\end{eqnarray*}
This completes the proof.
\end{pf}

\section{Hitting probabilities for the stochastic wave equation:
The~Gaussian case}
\label{s4}

In this section, we consider the solution to equation (\ref{1.1}),
which is the $d$-dimensional
Gaussian random field defined by
\begin{equation}
\label{4.1}
u(t,x) = \int_0^t\int_{\mathbb{R}^k}G(t-r,x-y) \sigma M(\mathrm{d}s,\mathrm{d}y) ,
\qquad
(t,x)\in[0,T]\times\mathbb{R}^k.
\end{equation}
Since $\sigma$ is invertible, we may assume, as in \cite{dkn07}, that
$\sigma$ is the identity matrix. Note that, in this case,
$u(t,x)=(u_1(t,x), \ldots, u_d(t,x))$ with
\[
u_i(t,x)=\int_0^t \int_{\mathbb{R}^k} G(t-r,x-y) M^i(\mathrm{d}s,\mathrm{d}y) ,\qquad
(t,x)\in[0,T]\times\mathbb{R}^k,
\]
$i=1,\ldots,d$, and, therefore,
the component processes
$(u_i(t,x),  (t,x)\in[0,T]\times\mathbb{R}^k)$, $i=1,\ldots,d$,
are i.i.d.



Most of the results of this section require the following hypothesis.

\begin{enumerate}[(H$_\beta$)]
\item[(H$_\beta$)]  The spectral covariance measure $\mu
$ is absolutely continuous with respect to Lebesgue measure on $\mathbb
{R}^k$ and
its density is given by
\[
f(\xi)=\Vert\xi\Vert^{-k+\beta} ,\qquad   \beta\in  \,]0,2\wedge k[.
\]
Equivalently, $\Gamma(\mathrm{d}x)= C(k,\beta)\Vert x\Vert^{-\beta}\,\mathrm{d}x$ (see
\cite{do}). Note that $(\mathrm{H}_\beta)$ implies (\ref{H}).
\end{enumerate}

In the sequel, we fix a strictly positive real number $t_0$.
We first aim for lower bounds on hitting probabilities. For this, we
intend to apply Theorem \ref{t3.1}. 
The required upper bound on the joint densities will be obtained by
combining Proposition \ref{p2.1} and the next two results.

\begin{prop}
\label{p4.2} Assume ${(\mathrm{H}_\beta)}$ and fix $M>0$. There then exist
positive constants $C_1$, $C_2$ such that for any
$(t,x), (s,y)\in[t_0,T]\times[-M,M]^k$,
\begin{eqnarray}
\label{4.3}
C_1 (|t-s| + \Vert x-y\Vert)^{2-\beta}&\le& E\bigl(
\Vert u(t,x)-u(s,y)\Vert^2\bigr)\nonumber\\[-8pt]\\[-8pt]
& \le& C_2 (|t-s| + \Vert x-y\Vert)^{2-\beta}.\nonumber
\end{eqnarray}
\end{prop}

\begin{pf}
The structure of this proof is similar to that
of Lemma 4.2 in \cite{dkn07}, but the methods for obtaining the
estimates differ substantially. Without loss of generality, we will
assume that $d=1$. Let $R(x)= E(u(t,x)u(t,0))$ with $t\ge t_0$. We then have
\[
E\bigl(\bigl(u(t,x)-u(t,y)\bigr)^2\bigr)= 2
\bigl(R(0)-R(x-y)\bigr).
\]
Following the steps of the proof of Remark 5.2 in \cite{dss} with the
dimension $k=3$ replaced by an arbitrary value of $k$
and, therefore, the Riesz kernel $\Vert\xi\Vert^{-(3-\beta)}$
replaced by $\Vert\xi\Vert^{-(k-\beta)}$,
we obtain
\begin{equation}
\label{4.4}
R(0)-R(x) \le C \Vert x\Vert^{2-\beta}.
\end{equation}
We next fix $y\in\mathbb{R}^k$ and consider increments in time. Let
$t_0\le
s<t\le T$. Using (\ref{e20}) and (\ref{fg}),
we have
\[
E\bigl(\bigl(u(t,y)-u(s,y)\bigr)^2\bigr) = S_1(s,t)+S_2 (s,t)
\]
with
\begin{eqnarray*}
S_1(s,t)&=&\int_0^s \mathrm{d}r \int_{\mathbb{R}^k} \frac{\mathrm{d}\xi}{\Vert\xi
\Vert
^{k-\beta}}\frac{|\sin((t-r)\Vert\xi\Vert)-\sin((s-r)\Vert\xi
\Vert)|^2}{\Vert\xi\Vert^2} ,\\
S_2(s,t)&=&\int_s^t \mathrm{d}r \int_{\mathbb{R}^k} \frac{\mathrm{d}\xi}{\Vert\xi
\Vert
^{k-\beta}}\frac{\sin^2((t-r)\Vert\xi\Vert)}{\Vert\xi\Vert^2}.
\end{eqnarray*}
With the changes of variables $r\to s-r$ and $\xi\to(t-s)\xi$, along
with the trigonometric formula $\sin x-\sin y= 2\sin\frac{x-y}{2}\cos
\frac{x+y}{2}$,
we obtain
\begin{eqnarray*}
S_1(s,t)&\le&4\int_0^s \mathrm{d}r \int_{\mathbb{R}^k}\frac{\mathrm{d}\xi}{\Vert\xi
\Vert
^{k-\beta+2}}\sin^2\biggl(\frac{(t-s)\Vert\xi\Vert}{2}\biggr)\\
&=& 4 \int_0^s \mathrm{d}r (t-s)^{2-\beta}\int_{\mathbb{R}^k} \frac
{\mathrm{d}v}{\Vert
v\Vert^{k-\beta+2}}\sin^2\biggl(\frac{\Vert v\Vert}{2}\biggr)\\
&\le& C |t-s|^{2-\beta}.
\end{eqnarray*}
For the term $S_2(s,t)$, we consider the changes of variables $r\to
t-r$ and then $\xi\to r\xi$, which easily yield
\begin{eqnarray*}
S_2(s,t)&\le&\int_0^{t-s} \mathrm{d}r\,  r^{2-\beta}\int_{\mathbb{R}^k} \frac
{\mathrm{d}v}{\Vert v\Vert^{k-\beta+2}} \sin^2\Vert v\Vert\\
&\le& C |t-s|^{3-\beta}.
\end{eqnarray*}
Hence, we have proven that
\begin{equation}
\label{4.5}
E\bigl(\bigl(u(t,y)-u(s,y)\bigr)^2\bigr) \le C |t-s|^{2-\beta}
\end{equation}
with a positive constant $C$ depending only on $T$.
With (\ref{4.4}) and (\ref{4.5}), we have established the upper bound
in (\ref{4.3}).

We now prove the lower bound in (\ref{4.3}) using several steps.

\begin{longlist}
\item[\textit{Step 1}.]
Assume that $s=t\ge t_0$ and $x\neq y$.
The arguments in the proof of Theorem 5.1(a) in \cite{dss} can be
trivially extended to any spatial dimension $k$.
Therefore, there is a positive constant $c_1$ such that for any $x, y
\in[-M,M]^k$,
\begin{equation}
\label{4.6}
E\bigl(\bigl(u(t,x)-u(t,y)\bigr)^2\bigr)\ge c_1 |x-y|^{2-\beta}.
\end{equation}

\item[\textit{Step 2}.]
We show that for arbitrary $x,y\in[-M,M]^k$ and $t_0\le s\le t\le T$,
\begin{equation}
\label{4.7}
E\bigl(\bigl(u(t,x)-u(s,y)\bigr)^2\bigr)\ge c |t-s|^{2-\beta}.
\end{equation}

Indeed, the left-hand side of this inequality is equal to
\[
R_1(s,t;x,y)+R_2(s,t;x,y)
\]
with
\begin{eqnarray*}
R_1(s,t;x,y)
& =&\int_0^s \mathrm{d}r \int_{\mathbb{R}^k} \frac{\mathrm{d}\xi}{\Vert\xi\Vert
^{k-\beta
}}\vert\mathcal{F}G(t-r,x-\cdot)(\xi
)-\mathcal{F}
G(s-r,y-\cdot)(\xi)\vert^2,\\
R_2(s,t;x,y)&=&\int_s^t \mathrm{d}r \int_{\mathbb{R}^k} \frac{\mathrm{d}\xi}{\Vert
\xi\Vert
^{k-\beta}}\vert\mathcal{F}G(t-r,x-\cdot)(\xi)\vert^2.
\end{eqnarray*}
Since $R_2(s,t;x,y)$ is positive, we can neglect its contribution. (We
note that
\[
R_2(s,t;x,y)\ge C |t-s|^{3-\beta}
\]
for some positive constant $C$. For $k=3$, this is shown in the proof
of Theorem 5.1 in \cite{dss} and it is easy to check that
the arguments go through to any dimension.)

By developing the integrand in $R_1(s,t;x,y)$, we find that
\begin{eqnarray*}
&&\Vert\xi\Vert^2\vert\mathcal{F}G(t-r,x-\cdot)(\xi
)-\mathcal{F}
G(s-r,y-\cdot)(\xi)\vert^2\\
&&\quad   = \big\vert\sin\bigl((t-r)\Vert\xi\Vert\bigr)- \mathrm{e}^{\mathrm{i}\xi\cdot
(y-x)}\sin\bigl((s-r)\Vert\xi\Vert\bigr)\big\vert^2\\
&&\quad  =\frac{1-\cos(2(t-r)\Vert\xi\Vert)}{2}+\frac{1-\cos
(2(s-r)\Vert\xi\Vert)}{2}\\
&&\qquad {}  -\cos\bigl(\xi\cdot(y-x)\bigr)\bigl[\cos\bigl((t-s)\Vert\xi\Vert
\bigr)-\cos\bigl((t+s-2r)\Vert\xi\Vert\bigr)\bigr].
\end{eqnarray*}
After integrating this last expression with respect to the variable
$r$, we obtain a positive quantity which is the sum of the following
three terms:
\begin{eqnarray*}
A_1&=&s\bigl[1-\cos\bigl((t-s)\Vert\xi\Vert\bigr) \cos\bigl(\xi\cdot
(y-x)\bigr)\bigr] ;\\
A_2&=&\frac{\sin((s+t)\Vert\xi\Vert)}{2\Vert\xi\Vert}\bigl(
\cos\bigl(\xi\cdot(y-x)\bigr)-\cos\bigl((t-s)\Vert\xi\Vert\bigr)\bigr)  ;\\
A_3&=&\frac{\sin(2(t-s)\Vert\xi\Vert)}{4\Vert\xi\Vert}-\frac
{\sin((t-s)\Vert\xi\Vert)}{2\Vert\xi\Vert} \cos\bigl(\xi\cdot
(y-x)\bigr).
\end{eqnarray*}
For the integration with respect to the variable $\xi$, we restrict
the domain to the set
\[
D_0= \bigl\{\xi\in\mathbb{R}^k\dvt \Vert\xi\Vert(t-s)\ge1,  \cos
\bigl((t-s)\Vert
\xi\Vert\bigr)\ge0\bigr\}.
\]
Note that on $D_0$, we have $A_1\ge0$. In fact,
\begin{eqnarray*}
A_1&=&s\bigl[1-\cos\bigl((t-s)\Vert\xi\Vert\bigr)+\cos\bigl((t-s)\Vert\xi\Vert
\bigr)\bigl(1-\cos\bigl(\xi\cdot(y-x)\bigr)\bigr)\bigr]\\
&\ge& s\bigl[1-\cos\bigl((t-s)\Vert\xi\Vert\bigr)\bigr].
\end{eqnarray*}
Moreover, 
\[
\vert A_2+A_3\vert\le\frac{2}{\Vert\xi\Vert}.
\]
Thus, with the change of variables $\xi\to(t-s)\xi$, we easily obtain
\begin{eqnarray*}
\int_{D_0} \frac{\mathrm{d}\xi}{\Vert\xi\Vert^{k-\beta+2}} A_1&\ge& s\int
_{D_0}\frac{\mathrm{d}\xi}{\Vert\xi\Vert^{k-\beta+2}}\bigl[1-\cos\bigl((t-s)\Vert
\xi\Vert\bigr)\bigr]\\
&=& s|t-s|^{2-\beta}\int_{\{\Vert w\Vert\ge1;  \cos(\Vert w\Vert
)\ge0\}} \frac{\mathrm{d}w}{\Vert w\Vert^{k-\beta+2}}\bigl(1-\cos(\Vert w\Vert
)\bigr)\\
&\ge& c_2 |t-s|^{2-\beta}.
\end{eqnarray*}
Similarly,
\[
\int_{\{\Vert\xi\Vert(t-s)\ge1\}} \frac{\mathrm{d}\xi}{\Vert\xi\Vert
^{k-\beta+2}} \vert A_2+A_3\vert\le c_3 |t-s|^{3-\beta}.
\]
Therefore, by the triangle inequality, we obtain
\[
R_1(s,t;x,y)\ge c_2|t-s|^{2-\beta}-c_3 |t-s|^{3-\beta} \ge\frac
{c_2}{2}|t-s|^{2-\beta}
\]
if $|t-s|\le\frac{c_2}{2c_3}$. This proves (\ref{4.7}) for small
values of $|t-s|$.

To extend the validity of (\ref{4.7}) to arbitrary values of $|t-s|$,
we note that
$R_1(s,t;x,y)$ is a continuous and positive function of its arguments
and, therefore, it is bounded below on
$\{(s,t;x,y)\in[t_0,T]^2\times[-M,M]^{2k}\dvt |t-s|\ge\varepsilon\}$
by some constant $c_\varepsilon$
for any $\varepsilon>0$. Hence, if $2T > |t-s|> \frac{c_2}{2c_3}$, we
also have
\[
R_1(s,t;x,y)\ge c |t-s|^{2-\beta}
\]
for some sufficiently small $c$.

\item[\textit{Step 3}.]
Suppose that $|t-s|\ge[\frac{c_1}{4C_2}]^{{1}/({2-\beta})}|x-y|$, where $c_1$ appears in (\ref{4.6}) and $C_2$ in
the right-hand side of (\ref{4.3}). By Step 2, we clearly have
\begin{eqnarray*}
E\bigl(\bigl(u(t,x)-u(s,y)\bigr)^2\bigr)&\ge& c |t-s|^{2-\beta}\\
&\ge& c\biggl(\frac{|t-s|}{2}+\frac{1}{2}\biggl(\frac
{c_1}{4C_2}\biggr)^{{1}/{(2-\beta)}}|x-y|\biggr)^{2-\beta}\\
&\ge& C_3( |t-s|+|x-y|)^{2-\beta}.
\end{eqnarray*}

\item[\textit{Step 4}.]
Suppose that $|t-s|\le[\frac{c_1}{4C_2}]^{{1}/{(2-\beta)}}|x-y|$. We then have
\begin{eqnarray*}
&&E\bigl(\bigl(u(t,x)-u(s,y)\bigr)^2\bigr)\\
&&\quad \ge\frac{1}{2} E
\bigl(\bigl(u(t,x)-u(t,y)\bigr)^2\bigr)-E\bigl(
\bigl(u(t,y)-u(s,y)\bigr)^2\bigr)\\
&&\quad \ge\frac{1}{2} c_1|x-y|^{2-\beta}-C_2|t-s|^{2-\beta}\\
&&\quad \ge\frac{c_1}{4}|x-y|^{2-\beta}\\
&&\quad \ge\frac{c_1}{4}\biggl(\frac{|x-y|}{2} + \frac{1}{2}\biggl[\frac
{4C_2}{c_1}\biggr]^{{1}/{(2-\beta)}}|t-s|\biggr)^{2-\beta}\\
&&\quad \ge C_4( |t-s|+|x-y|)^{2-\beta}.
\end{eqnarray*}
With this, the lower bound in (\ref{4.3}) is proved.\qed
\end{longlist}
\noqed\end{pf}

\begin{rem}
(a) As mentioned in the \hyperref[s1]{Introduction}, Proposition \ref{p4.2}, together
with Lemma \ref{l2.1}, establishes conditions (C1) and (C2) of \cite
{x} for the process
$U$.

(b) A consequence of the preceding proposition is that the sample paths
of (\ref{4.1}) are H\"older continuous, jointly in $(t,x)$, of
exponent $\gamma\in  \,]0, \frac{2-\beta}{2}[$, but they are not H\"
older continuous of exponent $\gamma> \frac{2-\beta}{2}$. We refer
the reader to \cite{dss} for a similar result on the solution to a
nonlinear stochastic wave equation in spatial dimension $k=3$.
\end{rem}

The next proposition is a further step toward proving that the process
$U$ satisfies the assumptions of Proposition \ref{p2.1}. We denote by
$\sigma^2_{t,x}$ the common variance of $u_i(t,x)$, $i=1,\ldots,d$.

\begin{prop}
\label{p4.1}
Assume that condition \textup{(\ref{H})} is satisfied. Fix $(t,x), (s,y)\in
[t_0,T]\times\mathbb{R}^k$. Then:
\begin{enumerate}[(ii)]
\item[(i)] $\sigma_{t,x}^2 \ge C(t_0\wedge t_0^3)>0$;
\item[(ii)] $|\sigma_{t,x}^2-\sigma_{s,y}^2| \le C|t-s|$.
\end{enumerate}
If, in addition, we assume that for $k^\prime<k$, all $k^\prime
$-dimensional submanifolds of $\mathbb{R}^k$ are sets with null $\mu
$-measure, then:
\begin{enumerate}[(iii)]
\item[(iii)] for any $(t,x)\neq(s,y)$ and $i=1,\ldots,d$,
\[
\operatorname{Corr}(u^i(t,x),u^i(s,y)) < 1.
\]
\end{enumerate}
\end{prop}

\begin{pf}
The variance of $u(t,x)$ is
\begin{equation}
\label{formvar}
\sigma_{t,x}^2= \int_0^t \mathrm{d}r \int_{\mathbb{R}^k} \mu(\mathrm{d}\xi)  \frac
{\sin
^2((t-r)\Vert\xi\Vert)}{\Vert\xi\Vert^2}
\end{equation}
and satisfies
\begin{equation}
\label{var}
C(t\wedge t^3)\le\sigma_{t,x}^2\le\bar C(t+t^3)
\end{equation}
(see, for instance, \cite{ss}, Lemma 8.6). This proves (i).

Assumption (\ref{H}) implies that
\[
\sup_{r\in[0,T]}\int_{\mathbb{R}^k} \mu(\mathrm{d}\xi) \frac{\sin
^2(r\Vert\xi
\Vert)}{\Vert\xi\Vert^2} \le C.
\]
Hence, assuming $t_0\le s<t\le T$, we obtain
\[
|\sigma_{t,x}^2-\sigma_{s,y}^2|=\int_s^t \mathrm{d}r \int_{\mathbb{R}^k}
\mu(\mathrm{d}\xi
) \frac{\sin^2(r\Vert\xi\Vert)}{\Vert\xi\Vert^2}\le C(t-s) ,
\]
which yields the conclusion (ii) of the proposition.

We now prove (iii) by checking that for any $(t,x)\neq(s,y)$ in
$[t_0,T]\times\mathbb{R}^k$,
\[
\sigma^2_{t,x}\sigma_{s,y}^2 - \sigma_{t,x;s,y}^2 >0 ,
\]
where $\sigma_{t,x;s,y}$ denotes the covariance of $u_i(t,x)$ and
$u_i(s,y)$ for any $i=1,\ldots,d$.

\begin{longlist}
\item[\textit{Case 1}:] $s<t$.
If $\sigma^2_{t,x}\sigma_{s,y}^2 - \sigma_{t,x;s,y}^2$ were equal to
zero, then the random variables
$u^i(t,x)$ and $u^i(s,y)$ would have correlation equal to $1$;
therefore, there would be $\lambda\in\mathbb{R}$
such that $u^i(t,x) = \lambda u^i(s,y)$ a.s.~and, in particular, we
would have
\[
E\bigl(\bigl(u^i(t,x)-\lambda u^i(s,y)\bigr)^2\bigr) = 0.
\]
The left-hand side of this equality is
\begin{eqnarray*}
&& \int_s^t \mathrm{d}r \int_{\mathbb{R}^k} \mu(\mathrm{d}\xi) \vert\mathcal
{F}G(t-r,x-\cdot
)(\xi)\vert^2\\
&&\quad {} + \int_0^s \mathrm{d}r \int_{\mathbb{R}^k} \mu(\mathrm{d}\xi) \vert
\mathcal{F}
G(s-r,x-\cdot)(\xi)-\lambda\mathcal{F}G(s-r,y-\cdot)(\xi)
\vert^2,
\end{eqnarray*}
which is bounded below, as in \eqref{var}, by $C((t-s)\wedge
(t-s)^3)$. This leads to a contradiction.

\item[\textit{Case 2}:] $s=t$, $x\ne y$.
We start, as in the preceding case, by assuming that $\sigma
^2_{t,x}\sigma_{t,y}^2 -\break  \sigma_{t,x;t,y}^2=0$,
hence
\[
E\bigl(\bigl(u^i(t,x)-\lambda u^i(t,y)\bigr)^2\bigr) = 0
\]
for some $\lambda\in\mathbb{R}$.
The left-hand side is equal to
\[
\int_0^t \mathrm{d}r \int_{\mathbb{R}^k} \mu(\mathrm{d}\xi) \vert \mathrm{e}^{\mathrm{i}\xi
\cdot
x}-\lambda \mathrm{e}^{\mathrm{i}\xi\cdot y}\vert^2 \vert\mathcal
{F}G(r,\cdot
)(\xi)\vert^2.
\]
If $\lambda=1$, then the integrand vanishes when $\cos[\xi\cdot
(x-y)]=1$ or $\sin(r\Vert\xi\Vert)=0$, which occurs on a $(k-1)$-dimensional
manifold of $\mathbb{R}^k$. Hence, by the assumption on $\mu$, we
reach a
contradiction.

If $\lambda\ne1$, then
\begin{eqnarray*}
&&\int_0^t \mathrm{d}r \int_{\mathbb{R}^k} \mu(\mathrm{d}\xi) \vert \mathrm{e}^{\mathrm{i}\xi
\cdot
x}-\lambda \mathrm{e}^{\mathrm{i}\xi\cdot y}\vert^2 \vert\mathcal
{F}G(r,\cdot
)(\xi)\vert^2\\
&&\quad  \ge\int_0^t \mathrm{d}r \int_{\mathbb{R}^k} \mu(\mathrm{d}\xi) (1-\lambda
)^2 \frac{\sin^2(r\Vert\xi\Vert)}{\Vert\xi\Vert^2}.
\end{eqnarray*}
This last integrand vanishes only when $\sin(r\Vert\xi\Vert)=0$.
Thus, we also get a contradiction in this case.
The proof of the proposition is now complete.\qed
\end{longlist}
\noqed\end{pf}

We can now obtain the required properties on densities, as follows.

\begin{prop}
\label{p4.3}
Assume $(\mathrm{H}_\beta)$. Fix $M, N>0$ and $(t,x), (s,y)\in
[t_0,T]\times[-M,M]^k$ with $(t,x) \neq(s,y)$.
\begin{enumerate}[(b)]
\item[(a)]
Let $p_{t,x;s,y}(\cdot,\cdot)$ denote the joint density of the random
vector $(u(t,x), u(s,y))$. We then have
\begin{equation}
\label{4.8}
p_{t,x;s,y}(z_1,z_2)\le\frac{C}{(|t-s| + |x-y|)^{{d(2-\beta)}/{2}}}
\exp\biggl(-\frac{c\Vert z_1-z_2\Vert^2}{(|t-s| + |x-y|
)^{2-\beta}}\biggr)
\end{equation}
for any $z_1,z_2\in[-N,N]^d$, where $C$ and $c$ are positive constants
not depending on $(t,x)$, $(s,y)$.
\item[(b)] Let $p_{t,x}$ denote the density of the random vector $u(t,x)$.
Then, for each $(t,x)\in[t_0,T]\times\mathbb{R}^k$ and $z\in[-N,N]^d$,
\begin{equation}
\label{4.9}
p_{t,x}(z)\ge C
\end{equation}
and
\begin{equation}
\label{4.90}
\sup_{z\in[-N,N]^d}\sup_{(t,x)\in[t_0,T]\times\mathbb{R}^k}
p_{t,x}(z)\le C.
\end{equation}
\end{enumerate}
\end{prop}

\begin{pf}
By Propositions \ref{p4.1} and \ref{p4.2}, we
see that the process $U$ satisfies the hypotheses of Proposition \ref
{p2.1} with $\eta=\frac{\beta}{2-\beta}$. Thus, we have statement (a).

The density $p_{t,x}$ is given by
\[
p_{t,x}(z) = \frac{1}{(2\curpi\sigma^2_{t,x})^{{d}/{2}}} \exp
\biggl(-\frac{\Vert z\Vert^2}{2\sigma^2_{t,x}}\biggr)
\]
with $\sigma_{t,x}^2$ as in (\ref{formvar}).
By (\ref{var}), we obtain both (\ref{4.9}) and (\ref{4.90}).
\end{pf}

The next theorem gives lower bounds on hitting probabilities.

\begin{teorem}
\label{t4.1}
Assume $(\mathrm{H}_\beta)$. Let $I$, $J$ be compact subsets of $[t_0,T]$
and $\mathbb{R}^k$, respectively, each with positive Lebesgue measure. Fix
$N>0$. Then:
\begin{enumerate}[(2)]
\item[(1)] there exists a positive constant
$c=c(I,J,N,\beta,k,d)$ such that for any Borel set $A\subset[-N,N]^d$,
\begin{equation}
\label{4.10}
P\{u(I\times J)\cap A\neq\varnothing\} \ge c  \operatorname{Cap}_{d-{2(k+1)}/{(2-\beta)}}(A) ;
\end{equation}
\item[(2)] for any $t\in I$, there exists a positive constant
$c=c(J,N,\beta,k,d,t)$ such that, for any Borel set $A\subset[-N,N]^d$,
\begin{equation}
\label{4.11}
P\bigl\{u(\{t\}\times J)\cap A\neq\varnothing\bigr\}\ge c  \operatorname{Cap}_{d-{2k}/{(2-\beta)}}(A) ;
\end{equation}
\item[(3)] for any $x\in J$, there exists a positive constant
$c=c(I,N,\beta,k,d,x)$ such that for any Borel set $A\subset[-N,N]^d$,
\begin{equation}
\label{4.12}
P\bigl\{u(I\times\{x\})\cap A\neq\varnothing\bigr\}\ge c  \operatorname{Cap}_{d-{2}/{(2-\beta)}}(A).
\end{equation}
\end{enumerate}
\end{teorem}

\begin{pf}
The three statements follow from Theorem \ref{t3.1} and
Proposition \ref{p4.3} applied, respectively, to the stochastic
process $U$, $U(t)=\{u(t,x), x\in\mathbb{R}^k\}$ with $t\in I$,
and $U(x)=\{u(t,x), t\in[t_0,T]\}$ with $x\in J$. Note that by (\ref
{4.8}) and (\ref{4.9}), the parameters $\gamma$ and $\alpha$
in Theorem \ref{t3.1} are $\gamma=\frac{d(2-\beta)}{2}$, $\alpha
=2-\beta$ and $m=k+1$, $m=k$, $m=1$, respectively.
\end{pf}

\begin{rem}
Since the probability of visiting translates of a compact set $A$
decreases to $0$ as the distance of this translated set to the origin
tends to infinity, it is not possible to replace $[-N,N]^d$ by $\mathbb{R}^d$
in the above theorem. In contrast, this \textit{will} be possible in the
upper bounds of the next theorem.
\end{rem}


\begin{teorem}
\label{t4.2}
Assume ${(\mathrm{H}_\beta)}$. Let $I$, $J$ be compact subsets of $[t_0,T]$
and $\mathbb{R}^k$, respectively, each with positive Lebesgue measure. Then:
\begin{enumerate}[(3)]
\item[(1)] there exists a positive constant
$c=c(I,J,\beta,k,d)$ such that for any Borel set $A\subset\mathbb{R}^d$,
\begin{equation}
\label{4.13}
P\{u(I\times J)\cap A\neq\varnothing\} \le c  {\mathcal{H}}_{d-{2(k+1)}/{(2-\beta)}}(A) ;
\end{equation}
\item[(2)] for any $t\in I$, there exists a positive constant
$c=c(J,\beta,k,d,t)$ such that for any Borel set $A\subset\mathbb{R}^d$,
\begin{equation}
\label{4.14}
P\bigl\{u(\{t\}\times J)\cap A\neq\varnothing\bigr\}\le c {\mathcal{H}}_{d-{2k}/{(2-\beta)}}(A) ;
\end{equation}
\item[(3)] for any $x\in J$, there exists a positive constant
$c=c(I,\beta,k,d,x)$ such that for any Borel set $A\subset\mathbb{R}^d$,
\begin{equation}
\label{4.15}
P\bigl\{u(I\times\{x\})\cap A\neq\varnothing\bigr\}\le c {\mathcal{H}}_{d-{2}/{(2-\beta)}}(A).
\end{equation}
\end{enumerate}
\end{teorem}

\begin{pf}
We first note that if we replace $d$ in the Hausdorff
dimensions of the bounds by any $\gamma\in  \,]0,d[$, then these
statements would be
a consequence of Theorem \ref{t3.2} applied, respectively, to the stochastic
processes $U$, $U(t)=\{u(t,x), x\in\mathbb{R}^k\}$ with $t\in I$,
and $U(x)=\{u(t,x), t\in[t_0,T]\}$ with $x\in J$. Indeed,
assumption (1) of Theorem \ref{t3.2}
is given in (\ref{4.90}). Moreover, since $U$ is a Gaussian process,
the right-hand side of (\ref{4.3}) yields the
validity of hypothesis (2) of Theorem~\ref{t3.2}, with $\delta=\frac
{2-\beta}{2}$.

The improvement to $\gamma=d$ is obtained by applying Theorem \ref
{t3.3} to each of the stochastic processes mentioned
before. Let us argue with the process $U$, for the sake of illustration.
From~(\ref{4.3}), we easily deduce that
\[
E\biggl[\exp\biggl\{\frac{\vert u_i(s,y)-u_i(t,x)\vert}{
(|s-t|+|x-y|)^{{(2-\beta)}/{2}}}\biggr\}\biggr]\le E
[\exp(cX)]=C ,
\]
where $X$ stands for a standard Normal random variable. Thus, when
$m=k+1$ and $\delta=\frac{2-\beta}{2}$, the left-hand side of (\ref
{m}) is bounded by a constant times
the square of the volume of $R_j^\varepsilon$, that is, $C\varepsilon
^{{4(k+1)}/{(2-\beta)}}$. Hence, the assumptions of Theorem \ref{t3.3}
are satisfied.

The proof of the theorem is complete.
\end{pf}

\section*{Acknowledgements}

 Research of Robert C. Dalang was supported in part by the
Swiss National Foundation for Scientific Research. Research of Marta
Sanz-Sol\'e was supported by the Grant MTM
2006-01351 from the Direcci\'on General de Investigaci\'on, Ministerio
de Ciencia e Innovaci\'on, Spain.

\printhistory

\end{document}